\documentclass{amsart}

\usepackage[mathcal]{euscript}
\usepackage{mathrsfs}
\usepackage{amsfonts,amssymb,amsbsy,amsmath}
\usepackage{graphicx,color}
\usepackage{epsfig}
\usepackage[all,cmtip]{xy}
\usepackage{verbatim}

\newcommand{\R}{\mathbb{R}}

\newcommand{\Z}{\mathbb{Z}}

\newcommand{\ba}{\begin{array}}
\newcommand{\ea}{\end{array}}

\newcommand{\GM}{\textup{Map}(T^2 \smallsetminus \textup{int}{D^2}, \partial)}
\newcommand{\M}{\textup{Map($T^2$)}}

\newcommand{\SL}{\textup{SL(2,\,$\Z$)}}
\newcommand{\GL}{\textup{GL(2,\,$\Z$)}} 
\theoremstyle{plain}
\newtheorem{theorem}{Theorem}
\newtheorem{lemma}[theorem]{Lemma}
\newtheorem{prop}[theorem]{Proposition}

\theoremstyle{remark}
\newtheorem{remark}{Remark}

\theoremstyle{definition}

\begin{document}
\title{A real open book not fillable by a real Lefschetz fibration}

\author{Fer\.{\i}t \"{O}zt\"{u}rk}
\address{Bo\u{g}az\.{\i}\c{c}\.{\i} \"{U}n\.{\i}vers\.{\i}tes\.{\i}, Department of Mathematics, TR-34342
  Bebek, \.Istanbul, Turkey}
\email{ferit.ozturk@boun.edu.tr}

\author{Nerm\.{\i}n Salepc\.{\i}}
\address{Institut Camille Jordan,
Universit\'e Lyon I,
43, Boulevard du 11 Novembre 1918
69622 Villeurbanne Cedex, France}
\email{salepci@math.univ-lyon1.fr}

\begin{abstract}
A real 3- or 4-manifold has by definition an orientation preserving smooth involution acting on it.
We consider  Lefschetz fibrations of 4-dimensional manifolds-with-boundary and open book decompositions on their boundary in the existence of a real structure.  
We prove that there is a real open book which cannot be filled by a real Lefschetz fibration, although it is filled by  non-real Lefschetz fibrations.
\end{abstract}

\maketitle

\section{Introduction and basic definitions}

There is a very close relation between open book decompositions and  Lefschetz fibrations over $D^2$, with fibers surfaces with boundary.   
It is known that any open book decomposition is filled by a Lefschetz fibration if its monodromy 
can be factorized  as a product of positive Dehn twists (work of Eliashberg, Giroux, Gompf, Loi-Piergallini, Akbulut-Ozbagci etc. See e.g.  \cite{lp}). Moreover, two fibrations filling the same open book are isomorphic if and only if  the two such factorizations  are Hurwitz equivalent \cite{kas}. 

In the present work we consider Lefschetz fibrations and open book decompositions on manifolds with a real structure.
A {\em real structure} on an oriented $4$-manifold (respectively $3$-manifold) with boundary (possibly empty) is defined 
to be an involution which is orientation preserving  and the fixed point set of which is  of  dimension $2$ (respectively $1$), if  it is not empty. Hence, if $M$ is the oriented boundary of an oriented  $4$-manifold $X$, a real 
structure on $X$ restricts to a real structure on $M$.
We call a manifold $X$ together with a real structure a {\em real manifold} and denote by $c_X$ the real structure on it.

Let $(L, \pi)$ be an open book decomposition of a real 3-manifold $(M, c_{M})$. Here $L\subset M$ is the binding 
and $\pi:M-L\rightarrow S^1$ is a fibration with fibers the page surfaces. We say that $(L, \pi)$ is a {\em real open book decomposition},  
if $\rho  \circ \pi=\pi \circ c_{M}|_{M-L}$ where $\rho: S^1 \to S^1$ is a reflection \cite{os2}. 

Let $p:E\to B$  be a genus-$g$ Lefschetz fibration of a 4-manifold $E$.
A  \emph{real Lefschetz fibration} is a Lefschetz fibration together with a pair of real structures $c_{E}:E\to E$ and $c_{B}: B\to B$ 
commuting with the fiber structure \cite{n2}.

In this article we answer a natural question that relates real open books and real Lefschetz fibrations.
We show that there is a Lefschetz fibration that does not admit a real structure, while
the canonical open book on its boundary is real. We also show that this real open book cannot be filled by any real Lefschetz
fibration with the same fiber topology and with arbitrary number of singular fibers (Theorem~\ref{dolmaz}).

\textbf{Acknowledgements.} The second author would like to thank  A.~Degtyarev for helpful discussions. 
The work in this article 
is supported by the Scientific and Technological Research Council of Turkey [TUBITAK-ARDEB-109T671].
The second author is supported by the European Community's Seventh Framework Programme ([FP7/2007-2013] [FP7/2007-2011]) under grant agreement no~[258204]. 

%%%%%%%%%%%%%%%%%%%%%%%%%%%%%%%%%%%REEL LEFSCHETZ LIFLENMLERI ILE ILISKI%%%%%%%%%%%%%%%%%%%%%%%%%%%%%%%%%%%%%%%%%%%%%%%%%%%%%%%%%%%%%%%%%%%%%%%%%%%%%%%%%%%%%%%%%%%%%%%%%%%%%%%%%%%%%%%%%%%%%%%%%%%%%%%%%%%%%%%%%%%%%%%%%%%%%

\section{Construction}

Let  $E$ be an oriented  smooth 4-manifold and  $B$ an oriented  smooth surface.
A \emph{genus-$g$ Lefschetz fibration} of  $E$ is a proper smooth projection $p:E\to B$ 
such that $p$ has only finitely many critical points %$y_{1},\ldots, y_{n}$ 
in int($E$) with pairwise distinct images %$b_{1},\ldots, b_{n}$ 
around which one can choose complex charts % $U_{i}\subset \C^2$ and $V_{i}\subset \C$ 
such that the projection takes the form $(z_{1},z_{2})\to z_{1}^2+z_{2}^2$. Moreover, the inverse image of a regular value is a closed oriented smooth surface of genus~$g$. It follows from the definition that if $\partial B\neq \varnothing$, then $\partial E=p^{-1}(\partial B)$ is a fiber bundle over $\partial B$.
The notion of Lefschetz fibration can be slightly generalized
to cover the case of fibers with boundary. Then $E$ turns into a
manifold with corners and its boundary, $\partial E$, becomes
naturally divided into two parts:  $p^{-1}(\partial B)$ % (called  the \emph{vertical boundary}) 
and $\partial F\times B$. %(called the \emph{horizontal boundary}). 
In this case, if,  in particular, $B=D^2$, then $\partial E$ admits a canonical open book decomposition with binding  $p^{-1}(0)$ and projection $p|_{\partial E}$.
%We call fibrations whose fibers are surfaces with boundary  \emph{Lefschetz fibrations with boundary}.  Fibrations with 
A  \emph{real Lefschetz fibration} is a fibration together with a pair of real structures $c_{E}:E\to E$ and $c_{B}: B\to B$ commuting with the fiber structure.  

It is known that around a singular fiber,  Lefschetz fibrations are  determined by the monodromy 
that is a single positive Dehn twist around a simple closed curve, the \emph{vanishing cycle}.  Algebraically, Lefschetz fibrations over $B=D^2$  can be encoded by the factorization $(t_{a_1}, \dots, t_{a_n})$ of the monodromy  $f=t_{a_n}\circ \ldots \circ t_{a_1}$ (because of the composition notation the order is reversed) along $\partial D^2$ into a product of positive Dehn twists considered up to  Hurwitz equivalence.  Namely, two such factorizations are called \emph{Hurwitz equivalent} if one can get  from one factorization to the other by a finite sequence of Hurwitz moves:
$$\begin{array}{lcr}(\ldots, t_{a_{i}}, t_{a_{i+1}},\ldots)  \to (\ldots,  t_{a_{i}}^{-1}\circ t_{a_{i+1}}\circ t_{a_{i}}, t_{a_{i}},\ldots),\\
(\ldots, t_{a_{i}}, t_{a_{i+1}},\ldots)  \to (\ldots,  t_{a_{i+1}}, t_{a_{i+1}}\circ t_{a_{i}}\circ t_{a_{i+1}}^{-1},\ldots);\end{array}$$
%they are called \emph{weakly Hurwitz equivalent} if they are equivalent via Hurwitz moves 
and possibly a global conjugation.

In what follows, we consider a genus-$1$ Lefschetz fibration $p: X\to D^2$  with exactly two singular fibers. Choose a base point $d\in \partial D^2$ and  consider a basis   $(\gamma_{1}, \gamma_{2})$ of $\pi_{1}(D^2\setminus\{\mbox{critical values}\}, d)$ obtained by  connecting  the base point  $d$ to the positively oriented simple loops, each surrounding the corresponding critical value once.  
Fix an identification of $F_{d}=p^{-1}(d)$  with $T^2=\R^2/\Z^2$.  Denote by $a$, the class  on $T^2$ of $(1,0)\in \R^2$ and by $b$, the class of $(0,1)$ %on $T^2$ 
and consider the curves $u, v$ represented respectively by $3a+5b$ and  $a$ on $T^2$.  We assume that the monodromy along $\gamma_{1}$
and $\gamma_{2}$ are given respectively by  the positive Dehn twists $t_{u}$ and $t_{v}$. 
In other words, the corresponding  singular fibers are obtained from $T^2$ by pinching the curves $u$ and $v$. 
The total monodromy of the fibration, thus, becomes the composition $f=t_{v}\circ t_{u}$.

Recall that $f\mapsto f_{*}$ defines an isomorphism from the mapping class group $\M$ (the identity component of the space of diffeomorphisms) of the torus to the group of automorphisms of $H_1(T^2,\Z)\cong\Z a\oplus\Z b$.  The latter is isomorphic to 
$$\SL =\left\{ [t_{a}]=\small{\left[\ba{rr}1 & 1 \\ 0 & 1\ea\right]}, [t_{b}]=\small{\left[\ba{rr} 1 & 0 \\ -1 & 1\ea\right]}: \begin{array}{ll} [t_a][t_b][ t_a]= [t_b][t_a][t_b]\\ ([t_a][t_b])^6=\textrm{id} \end{array}  \right\}$$
 where $[t_{\circ}]$ refers to the matrix representation of the automorphism $t_{\circ*}$. With respect to the above presentation,  $[f]=[t_{v}][t_{u}]$ is given by the matrix  $\small{\left[\ba{rr}-39 & 25 \\ -25 & 16\ea\right]}$.

\begin{prop}\label{reeldegil} $p:X\to D^2$ does not admit a real structure.
\end{prop}

\begin{figure}[h]
   \begin{center}
\includegraphics[scale=0.4]{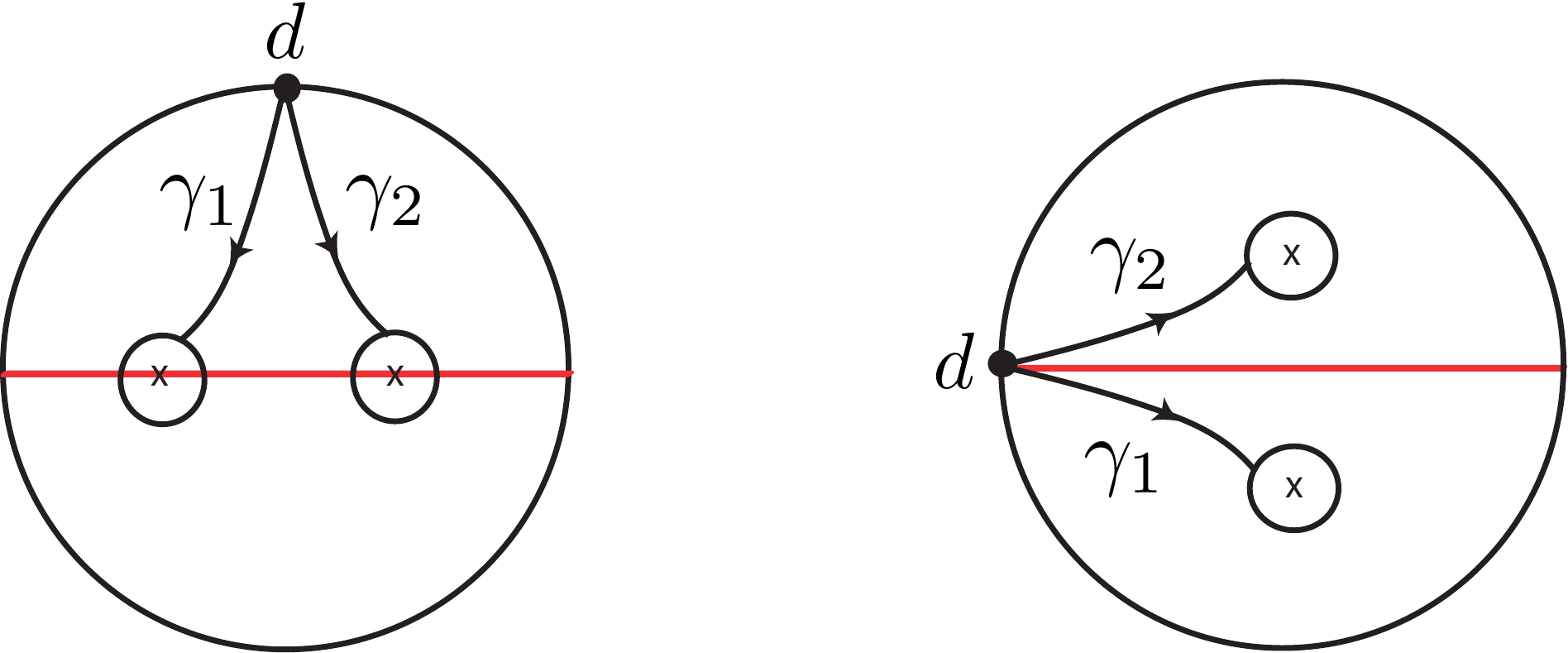}
\caption{Two possible positions of $(\gamma_{1},\gamma_{2})$ on base $D^2$ 
regarding the real structure reflection with respect to horizontal.}
       \label{gamalar}
      \end{center}
\end{figure}

\noindent {\it Proof:} Suppose that  $p: X\to D^2$ admits a real structure. That is to say, there exist real structures $c_{X}:X\to X$ and $c_{D^2}:D^2\to D^2$ such that $p\circ c_{X}=c_{D^2} \circ p$. Therefore, the critical points as well as their images are invariant under the action of the real structures $c_{X}$ and $c_{D^2}$, respectively, so either they are both real or they are interchanged by the real structures.  By definition of real Lefschetz fibrations, the decomposition $(t_{u},t_{v} )$ associated to $(\gamma_{1}, \gamma_{2})$ is Hurwitz equivalent to a decomposition $(t_{c_{X}(v)},t_{c_{X}(u)})$ associated to $(c_{D^2}(\gamma_{2}), c_{D^2}(\gamma_{1}))$.  
Although, the positions of $(\gamma_{1},\gamma_{2})$ can be arbitrary with respect to the real structure $c_{D^2}$, by conjugation by a power of $f$ (such a conjugation preserves %strong
 Hurwitz classes) and by moving, if necessary, the base point on the boundary, we can assume that we have only the two positions shown in Figure~\ref{gamalar}.      
%However, the set $\{\gamma_{1}, \gamma_{2}\)$ is probably not invariant under the action of the real structure. 
In the former case (the case shown on the left in Figure~\ref{gamalar}),  the fiber $F_{d}$ can be endowed with a real structure $c$ by pulling a real structure on a real fiber  between the two real singular fibers. Up to isotopy, we can assume that  vanishing cycles  $u$, $v$  are invariant under the action of the real structure $c$. 
Lemma~\ref{kesismelerencok2}  concerns the intersection number of invariant curves and prohibits this case, (since the intersection number  of $u$ and $v$  is  equal to 5). 

Now assume that critical points are exchanged by the real structure (this is the case depicted on the right in Figure~\ref{gamalar}). 
Let  $c: T^2\to T^2$ be  the real structure conjugate to the inherited real structure on the base fiber $F_{d}$.
Then, by the assumption we have $c(u)=v$; as a consequence, $u+v$ and $u-v$ are elements of  $\pm$ eigenspaces $H^c_+$ and, respectively, $H^c_{-}$ of $c_{*} : H_1(T^2,\Z)\to H_1(T^2,\Z)$. Thus,  the primitive forms of the classes of  $u+v$ and $u-v$ form a basis, respectively of $H^c_{+}$ and $H^c_{-}$. Moreover, for  any real structure $c$, the intersection number of the bases of $H_{+}^c$ and $H_{-}^c$ is at most 2 which is violated by our choice of $u$ and $v$.
\hfill  $\Box$ \\

\begin{lemma}\label{kesismelerencok2} The number of intersection points of two invariant curves on a real torus can only be 0, 1 or 2.
\end{lemma}

\noindent {\it Proof:}
The real structure restricted to an invariant curve defines an action on the curve. This action can be   
%For each real structure $c: T^2\to T^2$, if  $c(a)=a$ then $c$ induces an action on $a$ and  such an action can be
the identity, a reflection or an antipodal involution; let us call those curves  a real curve, a reflection curve and an antipodal curve respectively. 
Recall that up to equivariant diffeomorphisms, the real structure on a torus are distinguished by the number of real components that can be 0,1 or 2.
for each real structure we have

\begin{itemize}
 \item if $c$ has no real components, then there exist two $c$-equivariant isotopy classes of antipodal curves and no classes of other types; 
\item if $c$  has one real component, then $T^2$ contains a unique $c$-equivariant isotopy class of non-contractible reflection curves, a unique class of
antipodal curves, and a unique real curve;
\item if $c$ has two real components, then $T^2$ contains no $c$-equivariant isotopy class of antipodal curves, a unique class of
reflection curves, and two classes of real curves.
\end{itemize}

The result follows from the following observations: a reflection curve has two real points, so the intersection with a real curve must be at  0, 1 or 2 points. Moreover, a reflection curve must intersect an antipodal curve at 0 or 2 points, and real curves and antipodal curves are disjoint.  
\hfill  $\Box$ \\
%%%%%%%%%%%%%%%%%%%%%%%%%%%%%%%%%%%%%%%%
%$p:M\rightarrow S^1$ is a  surface bundle. We observe
\begin{prop}\label{reelkenar}
Let $M= \partial X$; then $M$ as a surface bundle is real.%, i.e. the monodromy diffeomorphism of $p$ is real.
\end{prop}

\noindent {\it Proof:}
By \cite[Proposition~4]{n1}, it is enough to check whether its monodromy is real, i.e.  it admits a decomposition into a product of two real structures. Real structures on $T^2$ correspond to involutive elements of determinant $-1$ in $\GL$. Here we have  $[f]=\small{\left[\ba{rr}-39 & 25 \\ -25 & 16\ea\right]}=\small{\left[\ba{rr}5 & -3 \\ 8 & -5\ea\right]\left[\ba{rr}-120 & 77 \\ -187 & 120\ea\right]}$, a product of two involutive elements of determinant $-1$. 
\hfill  $\Box$ \\
\begin{remark}
In \cite{n1}, an explicit classification of real elements of $\SL$ as well as a particular decomposition of a product of two real structures  for each conjugacy class is given. According to \cite[Proposition~4]{n1}, a hyperbolic element, i.e. an element with absolute value of trace greater than $2$, of $\SL$  is real if and only if its \emph{cutting period cycle} $[a_{1},\ldots, a_{2k}]$ is odd-bipalindromic. 
For hyperbolic matrices, the cutting period cycle together with the sign of the trace is a complete invariant of conjugacy classes and it can be recovered from the trace and the periodic tail of the continued fraction expansion of the slope of the eigenvectors. In our specific example the cutting period cycle of $f$ is found to be $[1313]$  which is odd-bipalindromic, i.e. up to cyclic ordering, it has two palindromic pieces $1$, $313$ of odd length. \end{remark}

So far we have constructed a $T^2$-fibered Lefschetz fibration  cannot be real while the $T^2$ bundle at its boundary is real.
Now we show that a similar construction can be cooked up to obtain a real open book at the boundary. 

Denote by $\breve{p}: \breve{X} \to D^2$  the Lefschetz fibration with boundary obtained from $p:X\to D^2$ above by taking out  a neighborhood of a section. Note that such a section always exists for genus~1 Lefschetz fibrations \cite{Mo}.
%$\nu(q)\times  D^2$  where $\nu(q)$ is an open neighborhood of a regular real point $q\in F_{d}\setminus \{u, v\}$. (Away from the singular points the fibration is trivial.) % so it is meaningful to cut out $\nu(q)\times  D^2$.) 
The regular fibers of $\breve{p}$ are tori with one boundary component. % on which the real structure acts as reflection (since $q$ is a real point).  
The fibration $\breve{p}$ has no real structure, since otherwise $p$ would have one. % has a real structure if and only if  has one.  Thus, by Proposition~\ref{reeldegil}, $\breve{p}$ does not admit a real structure. 
However,  we have:

\begin{prop} The canonical open book of $\breve {M}=\partial\breve {X}$ admits a real structure. 
\end{prop}

\noindent {\it Proof:}   As a consequence of  Lemma~1 in \cite{os2} it is enough to show that the monodromy of $\breve{p}|_{\breve{M}}$ is real.
Note that the monodromy  of $\breve{p}|_{\breve{M}}$ is an element of  $\GM$, the group of relative isotopy classes of orientation preserving diffeomorphisms which are identity on the boundary. 

The crucial observation is that  being real in $\M$ is equivalent to being real in $\GM$.
%$$\GM=\left\{[t_\bold{a}], [t_\bold{b}]: [t_\bold{a}][t_\bold{b}][t_\bold{a}]= [t_\bold{b}][t_\bold{a}][t_\bold{b}]\right\}.$$ 
Namely, there is a short exact sequence 
$$0 \to \langle t_{\partial}\rangle \to \GM  \to \M \to 0$$
given by the central extension of $\M$, where $t_{\partial}$ is the Dehn twist along the boundary component.  
Obviously, images of real elements of $\GM$ are real.
For the converse, suppose im$[f]=[c][c']$ for some $[f]\in \GM$ and real structures $[c],[c'] \in \GL$.  If necessary, by replacing $c$ by $c\circ g$ and $c'$ by $g^{-1}\circ c'$ for some diffeomorphism $g$, we can assume that  $c$ and hence $c'$ leaves $\partial T^2$ invariant. We have $t_{\partial}=c' t_{\partial}^{-1}c'$. Thus  $t_{\partial}^{-1}c$  is a real structure which we denote by $c''$. 
As $[f]=[c][c'][t_{\partial}]^k$, then for some
integer $k$,   $[f]=[c] [c'] [c' \circ  c'']^k=[c][c'' \circ (c'\circ c'')^{k-1}]$. Being conjugate to a real structure, $c''\circ (c'\circ c'')^{k-1}$ is a real structure.
Thus, $[f]$ is real. \hfill  $\Box$ \\

With  similar hands-on approach as the one above, we prove the following further result.  

\begin{prop}\label{baskaikiliyok}
Up to isomorphism preserving the identification $T^2\to F_{d}$,  
there are exactly two Lefschetz fibrations with exactly two singular fibers
 filling the  canonical open book on $\breve {M}$; moreover, neither of the fibrations admit a real structure.
%The  canonical open book on $\breve {M}$ cannot be filled by a real Lefschetz fibration over $D^2$ with the same fiber topology and with exactly two singular fibers.
\end{prop}

\noindent {\it Proof:}
Let $u',v'$ be two simple curves on $T^2$ such that $t_{v'} \circ t_{u'}=f$. We have $[f]=[t_{v}][t_{u}]=\small{\left[\ba{rr}-39 & 25 \\ -25 & 16\ea\right]}$. Suppose $v'$ is the curve defined by  $\alpha a +\beta b $, %$\alpha\bold{a}+ \beta \bold{b}$; 
then  a simple calculation gives $[t_{v'}]= \small{ \left[\ba{rr} 1-\alpha\beta & \alpha^2 \\ -\beta^2 & 1+\alpha\beta \ea\right]}$.
All Dehn twists are represented by matrices whose traces have absolute value equal to $2$. 
Therefore, from $[t_{u'}] =[t_{v'}^{-1} f]$, we get the identity $25\alpha^2+25\beta^2-55\alpha\beta = 25$. 
This quadratic Diophantine equation has solutions $(\pm 1,0)$, $(0,\pm 1)$ and (the transpose of) 
these vectors multiplied from left by the powers of the matrix $S=\small{\left[\ba{rr}-3 & 5 \\ -5 & 8\ea\right]}$ (see e.g 
the step-by-step computation of the online application \cite{dar}). 
Note that $u$ and $v$ are in the solution set. Note also that $u', v'$,  the curves with class $b$ and $5a+8b$ respectively,
are in the solution set too.
The pairs $(t_{u}, t_{v})$ and $(t_{u'},t_{v'})$ are not Hurwitz equivalent. Indeed in that case 
there would be an invertible matrix  $K$ with  determinant 1, such that  either 
$$K^{-1}[t^{-1}_{u}] [t_{v}][t_{u}] K=[t_{u'}] \mbox{ and } K^{-1} [t_{u}] K=[t_{v'}]$$
or 
$$K^{-1}[t^{-1}_{v}] [t_{u}][t_{v}] K=[t_{u'}]  \mbox{ and } K^{-1} [t_{v}] K=[t_{v'}].$$
However, by straightforward calculation one can show that there exists  no such $K$.

Note also that the matrices $S$ and $[f]$ commute;  hence  for any pair $(x, y)$ such that $f=t_{y}\circ t_{x}$, 
the factorization $(S^k x,S^k y)$ is Hurwitz equivalent (with a fixed identification) to the factorization  $(x,y)$ for every $k\in \Z$.

As a consequence, we have two Hurwitz equivalence classes given by the pairs 
$(t_{u},t_{v})$ and $(t_{u'}, t_{v'})$.
To finish, note that the number of intersection points of $u'$ and $v'$ is 5, so the proof of  Proposition~\ref{reeldegil} applies to $(u', v')$ to show that  the fibrations defined by the decomposition  $(t_{u'}, t_{v'})$ is not real.
\hfill  $\Box$ \\

%The above results combined with a delicate observation that has been pointed out to us by A.~Degtyarev prove the following

\begin{theorem}
\label{dolmaz}
The canonical open book on $\breve {M}$ cannot be filled by any real Lefschetz fibration with the same fiber topology and with arbitrary
number of singular fibers.
\end{theorem}

\noindent {\it Proof:}
We will show that $\breve{M}$ cannot be filled by a (real or non-real) Lefschetz fibration with the number of singular fibers 
different from 2. % by showing that its monodromy cannot admit arbitrary decomposition into a product of positive Dehn twists.
Recall that $$\GM=\left\{[t_a], [t_b]: [t_a][t_b][ t_a]= [t_b][t_a][t_b] \right\}.$$ Therefore, the homomorphism 
$$\begin{array}{lrcl}
deg : &\GM&\to&\Z\\
&t_{a},t_{b}&\mapsto&1
\end{array}$$ is well-defined; hence the number
of Dehn twists that may constitute a given element $h$ equals $deg(h)$. In our case,  $deg(f)=2$, so it cannot  admit a factorization as a product of an arbitrary number of positive Dehn twists other than 2.  Moreover, all possible factorizations into a product of two Dehn twists are shown above to be non-real.
 \hfill  $\Box$ \\

\begin{remark} 
Elements admitting a factorization into a product of two Dehn twists will be studied in \cite{an}. 
There the classification of such factorizations up to Hurwitz equivalence are presented and their relation to real structures are also elaborated.
\end{remark}

\end{document}